\documentstyle[11pt]{article}
  \newlength{\standardunitlength}
\newtheorem{cor}{Corollary} \newtheorem{lemma}{Lemma}
\newtheorem{theorem}{Theorem} \newtheorem{prop}{Proposition}
\newenvironment{proof}{\noindent {\sc Proof:}}{$\Box$ \vspace{2 ex}}

\begin{document}

\begin{center} {GL(n,q) and Increasing Subsequences in Nonuniform Random Permutations} \end{center}

\begin{center}
By Jason Fulman
\end{center}

\begin{center}
University of Pittsburgh Math Department
\end{center}

\begin{center}
301 Thackeray Hall
\end{center}

\begin{center}
Pittsburgh, PA 15260
\end{center}

\begin{center}
fulman@math.pitt.edu
\end{center}

\begin{center}
Version 1: September 12, 2001
\end{center}

\begin{center}
Version 2: September 23, 2001
\end{center}

\newpage

\begin{abstract} Connections between longest increasing subsequences in random permutations and eigenvalues of random matrices with complex entries have been intensely studied. This note applies properties of random elements of the {\it finite} general linear group to obtain results about the longest increasing and decreasing subsequences in non-uniform random permutations. \end{abstract}

\section{Introduction} \label{introduction}

	In recent years there has been serious interest in the
	relationship between increasing subsequences of random
	permutations and eigenvalues of random complex matrices from
	various ensembles. It is beyond the scope of this paper to
	survey the subject, but the connections are fascinating and
	relate to Painleve functions, Riemann surfaces, solitaire,
	interacting particle systems, point processes, quantum
	mechanics, Riemann-Hilbert problems, and more. Recent surveys
	include \cite{AD} and \cite{De}.

	The purpose of this note is to give first relationships between ``eigenvalues'' of elements of finite classical groups and longest increasing subsequences. Section \ref{matrix} recalls a probability measure $P_{n,q}$ on partitions of size $n$, explaining its group theoretic meaning. There is a simple formula for the distribution of the number of parts of a partition chosen from $P_{n,q}$. Using connections with the Rogers-Selberg identity, Section \ref{matrix} derives results on the distribution of the largest part of a partition chosen from $P_{n,q}$. Section \ref{matrix} closes by proving a combinatorially interesting monotonicity result.

	Section \ref{subsequences} recalls a measure $Q_{n,q}$ on
	partitions of size $n$ and explains its relationship with
	increasing and decreasing subsequences in non-uniform random permutations and with
	unipotent representations of the finite general linear
	groups. The measure $Q_{n,q}$ is a natural $q$-analog of the
	Plancherel measure of the symmetric group. Then it is proved
	that although $P_{n,q}$ and $Q_{n,q}$ are different, they are
	sufficiently similar that information about $P_{n,q}$ can be
	used to deduce information about $Q_{n,q}$. This gives results
	about the first row and first column under the measure $Q_{n,q}$, and hence
	about the longest increasing and decreasing subsequence of non-uniform
	permutations. We remark that as $q \rightarrow \infty$ the
	measures $P_{n,q}$ and $Q_{n,q}$ both converge to the point
	mass on the one row partition of size $n$. This behavior is
	qualitatively different from other models such as the usual
	Plancherel measure on the symmetric group. Throughout the
	paper we assume that $q \geq 2$ so that $P_{n,q}$ and
	$Q_{n,q}$ are close enough to be usefully compared.

	The distribution of the first row or column under the measure $Q_{n,q}$
	could be studied via Toeplitz determinants \cite{BDJ},\cite{TW} or by the point
	process approach of \cite{BOO}. The approach here yields
	different insights than these approaches would and gives explicit bounds for all $n$. It
	also avoids the issue of having to derandomize the variable
	$n$ which occurs in these other approaches. In any case, our
	purpose here is to illustrate connections with finite group
	theory.

\section{Notation and Lemmas} \label{notation}

	To begin we describe some standard notation about partitions which will be used throughout the paper. Let $\lambda$ be a partition of some non-negative integer $|\lambda|$ into parts
$\lambda_1 \geq \lambda_2 \geq \cdots$. Let $m_i(\lambda)$ be the
number of parts of $\lambda$ of size $i$, and let $\lambda'$ be the
transpose of $\lambda$ in the sense that $\lambda_i' =
m_i(\lambda) + m_{i+1}(\lambda) + \cdots$. It is also useful to
define the diagram associated to $\lambda$ as the set of points $(i,j)
\in Z^2$ such that $1 \leq j \leq \lambda_i$. We use the convention
that the row index $i$ increases as one goes downward and the column
index $j$ increases as one goes across. So the diagram of the
partition $(4331)$ is:

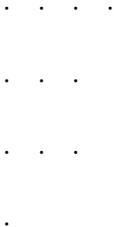
\[ \begin{array}{c c c c c}
                . & . & . & . &  \\
                . & . & . &  &    \\
                . & . & . &  &    \\
                . & & & &
          \end{array}  \] The hook length of a dot $s$ in $\lambda$ is defined as $a(s)+l(s)+1$ where $a(s)$ is the number of dots in the same row as $s$ to the right of $s$ and $l(s)$ is the number of dots in the same column of $s$ south of $s$.

	Throughout the paper we use the notation from $q$-series that $(x)_n=(1-x)(1-x/q)(1-x/q^2) \cdots (1-x/q^{n-1})$. We also use the following elementary lemmas.

\begin{lemma} \label{neumann1} (\cite{NP}) If $q \geq 2, d \geq 1$ then \[ (1-1/q)^2 \leq \prod_{i=1}^d (1-1/q^i) \leq 1-1/q .\] \end{lemma} 

	In fact \cite{NP} shows that for $q \geq 2$, $1-\frac{1}{q}-\frac{1}{q^2} \leq \prod_{i=1}^d (1-1/q^i)$. This strengthening would improve some of the bounds in this paper but we content ourselves with the bound from Lemma \ref{neumann1}.

\begin{lemma} \label{euler} (Euler) \[ \prod_{i \geq 1} (\frac{1}{1-u/q^i}) = 1+\sum_{n \geq 1} \frac{u^n}{q^n (1/q)_n}.\] \end{lemma}   

\section{The Measure $P_{n,q}$ on Partitions} \label{matrix}

	Recall that for the unitary group with complex entries $U(n,C)$, the set of eigenvalues of an element exactly parameterizes its conjugacy class. Hence it is natural to study conjugacy classes of a random element of $GL(n,q)$. A matrix $\alpha \in GL(n,q)$ uniquely decomposes the underlying vector space $V$ as a direct sum of subspaces $V_{\phi}$ where

\begin{enumerate}

\item $\phi$ is a monic irreducible polynomial with coefficients in the finite field $F_q$.

\item The characteristic polynomial of $\alpha$ restricted to $V_{\phi}$ is a power of $\phi$.

\item The characteristic polynomials of $\alpha$ restricted to distinct summands $V_{\phi_1}$ and $V_{\phi_2}$ are coprime.

\end{enumerate} Recall that a subspace $W$ invariant under $\alpha$ is called cyclic if it contains a vector $w$ such that $W$ is generated by $\{\alpha^i w, i \geq 0\}$. Each $V_{\phi}$ decomposes as a sum of a cyclic subspaces. Although this decomposition of $V_{\phi}$ need not be unique, the dimensions of the cyclic subspaces in the decomposition are uniquely determined and define a partition $\lambda_{\phi}(\alpha)$ where the parts of the partitions are the dimensions of the cyclic subspaces in the decomposition of $V_{\phi}$, each divided by the degree of $\phi$. Thus to each element $\alpha$ of $GL(n,q)$ is associated an infinite collection of partitions $\lambda_{\phi}(\alpha)$ and this data determines the conjugacy class of $\alpha$ \cite{H}. Note that one has the conditions that $\lambda_z$ is empty (since $\alpha$ is invertible) and that $\sum_{\phi} deg(\phi) |\lambda_{\phi}|=n$. Picking $\alpha$ uniformly at random in $GL(n,q)$ makes the $\lambda_{\phi}$ random variables.

	As $n \rightarrow \infty$, the random variables
	$\lambda_{\phi}$ become independent. Furthermore the law of
	$\lambda_{\phi}$ depends on $\phi$ only through its degree and
	in fact one can study $\lambda_{z-1}$ without loss of
	generality. Thus one has a very natural probability measure on
	the set of all partitions of all natural numbers. Further
	discussion of this measure can be found in the survey
	\cite{F2}. For our purposes we need the formula which says
	that the chance that this limit measure (which we denote
	$\tilde{P}_q$) yields $\lambda$ is \[ \prod_{i=1}^{\infty}
	(1-\frac{1}{q^i}) \frac{1}{\prod_{j \geq 1} q^{(\lambda_j')^2}
	(\frac{1}{q})_{m_j(\lambda)}}.\]

	The measure $P_{n,q}$ in this paper is given by renormalizing
	$\tilde{P}_q$ to live on partitions of size
	$n$. (This turns out to be equivalent to studying the random
	partition $\lambda_{z-1}$ for a uniformly chosen unipotent
	element of $GL(n,q)$. We note that it is {\it not} the same as
	looking at $\lambda_{z-1}$ for a uniformly chosen element
	$\alpha \in GL(n,q)$, since $\lambda_{z-1}(\alpha)$ could have
	size less than $n$).

\begin{prop} \label{Pform}
\[ P_{n,q}(\lambda) =  \frac{q^n (\frac{1}{q})_n}{\prod_j q^{(\lambda_j')^2} (\frac{1}{q})_{m_j(\lambda)}}.\]	
\end{prop}

\begin{proof} As the proof of Lemma 4 in Section 3.1 of the survey \cite{F2} explains, 

\[ \sum_{\lambda: |\lambda|=n} \frac{1}{\prod_j q^{(\lambda_j')^2} (\frac{1}{q})_{m_j(\lambda)}} = \frac{1}{q^n (\frac{1}{q})_n}.\] The result follows. \end{proof}

	Before discussing the distribution of the first row under $P_{n,q}$ we recall the distribution of the first column, for which there is a remarkably simple formula. (One can show that the size of first column of $\lambda_{z-1}(\alpha)$ is the dimension of the fixed space of $\alpha$).

\begin{theorem} \label{firstcolumn} (\cite{F1}) 

\begin{enumerate}

\item The $\tilde{P}_q$ probability that $\lambda_1'=k$ is $\prod_{i=1}^{\infty} (1-\frac{1}{q^i}) \frac{(1/q)^{k^2}}{(1/q)_k^2}$.

\item The $P_{n,q}$ probability that $\lambda_1'=k$ is $\frac{(1/q)_n (1/q)_{n-1}}{q^{k^2-k} (1/q)_k (1/q)_{k-1} (1/q)_{n-k}}$.

\end{enumerate}
\end{theorem}

	The remainder of this section studies the distribution of the first row under the measure $P_{n,q}$. We remark in passing that this statistic is interesting since the first row of $\lambda_{z-1}$ for a unipotent matrix determines the order of the matrix.	

	Let $P_{n,q}^r$ be the probability that the first row of a partition chosen from the measure $P_{n,q}$ has length strictly less than $r$. Proposition \ref{expand} gives an expansion for $P_{n,q}^r$.

\begin{prop} \label{expand}

\[ P_{n,q}^r = q^n (1/q)_n * \ Coeff. \ u^n \ in \ \prod_{i \geq 1} \frac{1}{1-u/q^i} \sum_{m=0}^{\infty} \frac{(-1)^m (1-u/q^{2m}) u^{rm} (u/q)_{m-1}}{q^{rm^2+{m \choose 2}} (1/q)_m} .\]

\end{prop}   

\begin{proof} Clearly 

\[ \sum_{\lambda: |\lambda|=n \atop \lambda_1<r} P_{n,q}(\lambda) = q^n (1/q)_n * \sum_{\lambda: |\lambda|=n \atop \lambda_1<r} \frac{1}{\prod_j q^{(\lambda_j')^2} (\frac{1}{q})_{m_j(\lambda)}}.\] Corollary 2 in Section 3.4 of \cite{F3} shows that

\[ \sum_{\lambda: \lambda_1<r} \prod_{i=1}^{\infty} (1-u/q^i) \frac{u^{|\lambda|}}{\prod_{j \geq 1} q^{(\lambda_j')^2} (\frac{1}{q})_{m_j(\lambda)}}\\
=  \sum_{m=0}^{\infty} \frac{(-1)^m (1-u/q^{2m}) u^{rm} (u/q)_{m-1}}{q^{rm^2+{m \choose 2}} (1/q)_m}.\]

\end{proof} 

	Now we prove the main result of this section. The assumption that $r \leq n-1$ is for convenience; it is simple to derive closed expressions for $P^r_{n,q}$ with $r \geq n$.

\begin{theorem} \label{Pbound} For $q \geq 2$ and $r \leq n-1$,

\begin{enumerate}
\item $P_{n,q}^r \leq  (\frac{1}{1-1/q})^2 [\frac{1}{q^{2n-2r+2}}+ \frac{1}{(1-1/q^{2n+1})}\frac{1}{q^{n+1}}]$

\item $P_{n,q}^r \geq \frac{1}{q^{2n-2r+2}}-\frac{1}{q^{n+1}(1-1/q)} - \frac{1}{q^{2n-2r+3}(1-1/q)^2} - \frac{1}{q^{3n-3r+4}(1-1/q)^3} - \frac{1}{(1-1/q)^2} \frac{1}{q^{2n+3} (1-1/q^{2n+3})}.$
\end{enumerate}
\end{theorem}

\begin{proof} First we prove the upper bound. Since $n \geq r+1$, Proposition \ref{expand} implies that

\begin{eqnarray*}
P_{n,q}^r & = & 1 + [ \frac{(-1) (1/q)_n}{(1/q)_{n-r} (1-1/q)} + \frac{(1/q)_n}{q (1/q)_{n-(r+1)} (1-1/q)} ]\\
& & +  [ q^n (1/q)_n \sum_{m \geq 2} Coeff. \ u^n \ in \ \prod_{i \geq 1} \frac{1}{1-u/q^i} \frac{(-1)^m (1-u/q^{2m}) u^{rm} (u/q)_{m-1}}{q^{rm^2+{m \choose 2}} (1/q)_m}].
\end{eqnarray*} Here the 1 comes from the $m=0$ term and the first term in
	square brackets comes from $m=1$. Consider the contribution
	from the $m=0$ and $m=1$ terms. It is equal to

\begin{eqnarray*}
& & 1 - (1-1/q^{n-r+1})(1-1/q^{n-r+2}) \cdots (1-1/q^n) (1+1/q^{n-r+1}+1/q^{n-r+2} +1/q^{n-r+3} + \cdots)\\
& \leq & 1 - (1-1/q^{n-r+1})(1-1/q^{n-r+2}) \cdots (1-1/q^n) (1+1/q^{n-r+1}+1/q^{n-r+2}+ \cdots + 1/q^{n})\\
& \leq & 1 - (1-1/q^{n-r+1}-1/q^{n-r+2}- \cdots - 1/q^{n}) (1+1/q^{n-r+1}+1/q^{n-r+2}+ \cdots + 1/q^{n})\\
& = & (1/q^{n-r+1}+ \cdots + 1/q^n)^2\\
& \leq & (\frac{1}{1-1/q})^2 \frac{1}{q^{2n-2r+2}}.
\end{eqnarray*} The second inequality used the fact that $(1-x_1) \cdots (1-x_r) \geq 1- (x_1+\cdots+x_r)$ if $0 \leq x_1,\cdots,x_r \leq 1$.

	To upper bound the second term in square brackets, observe that for $m \geq 2$ 

\[ \frac{q^n (1/q)_n (-1)^m}{q^{rm^2+{m \choose 2}} (1/q)_m} Coeff. \ u^{n-rm} \ in \ (1-u/q^{2m}) (u/q)_{m-1} \prod_{i \geq 1} \frac{1}{1-u/q^i} \] is positive only when $m$ is even, in which case it is less than

\begin{eqnarray*}
& &  \frac{q^n (1/q)_n}{q^{rm^2+{m \choose 2}} (1/q)_m} Coeff. \ u^{n-rm} \ in \ \prod_{i \geq m} \frac{1}{1-u/q^i}\\
& = & \frac{q^n (1/q)_n}{q^{rm^2+{m \choose 2}} (1/q)_m q^{m(n-rm)} (1/q)_{n-rm}}\\
& \leq & \frac{q^n}{q^{rm^2+{m \choose 2}+m(n-rm)}(1/q)_{n-rm}}\\
& \leq & \frac{1}{(1-1/q)^2} \frac{1}{q^{{m \choose 2} + n(m-1)}}.
\end{eqnarray*} The first inequality is true since $n \geq m$ and the second inequality is Lemma \ref{neumann1}.

	Thus the second term in square brackets is at most

\begin{eqnarray*}
\frac{1}{(1-1/q)^2} \sum_{m \geq 2 \atop m \ even} \frac{1}{q^{{m \choose 2} + n(m-1)}} & = & \frac{q^n}{(1-1/q)^2} \sum_{m \geq 1} \frac{1}{q^{2m^2-m+2mn}}\\
& \leq & \frac{q^n}{(1-1/q)^2} \sum_{m \geq 1} \frac{1}{q^{m(2n+1)}}\\
& = & (\frac{1}{1-1/q})^2 \frac{1}{q^{n+1}(1-1/q^{2n+1})}.
\end{eqnarray*}
  
	To lower bound $P_{n,q}^r$, we begin by examining the first term in square brackets. It is

\begin{eqnarray*}
& & 1-(1-1/q^{n-r+1}) \cdots (1-1/q^n)(1+1/q^{n-r+1}+1/q^{n-r+2}+\cdots+1/q^n+\cdots)\\
& \geq & 1-(1-(1/q^{n-r+1}+\cdots+1/q^n)+(\frac{1}{q^{n-r+1}} \frac{1}{q^{n-r+2}}+\frac{1}{q^{n-r+1}} \frac{1}{q^{n-r+3}}+\cdots+\frac{1}{q^{n-1}}\frac{1}{q^n})))\\
& & \cdot (1+1/q^{n-r+1}+1/q^{n-r+2}+\cdots+1/q^n+\cdots)\\
& \geq & 1-(1-(1/q^{n-r+1}+\cdots+1/q^n)+\frac{1}{q^{2n-2r+3}(1-1/q)^2})\\
& & \cdot (1+1/q^{n-r+1}+1/q^{n-r+2}+\cdots+1/q^n+\cdots)\\
& \geq & (1/q^{n-r+1}+\cdots+1/q^n)^2 -\frac{1}{q^{n+1}(1-1/q)} - \frac{1}{q^{2n-2r+3}(1-1/q)^2} - \frac{1}{q^{3n-3r+4}(1-1/q)^3}\\
& \geq & \frac{1}{q^{2n-2r+2}}-\frac{1}{q^{n+1}(1-1/q)} - \frac{1}{q^{2n-2r+3}(1-1/q)^2} - \frac{1}{q^{3n-3r+4}(1-1/q)^3}.
\end{eqnarray*}

	Note that the first inequality used the fact that $(1-x_1)\cdots(1-x_r) \leq 1-(x_1+\cdots+x_r)+(x_1x_2 + x_1x_3 + \cdots + x_{r-1}x_r)$ for $0 \leq x_1,\cdots,x_r \leq 1$. 

	Next we consider the second term in square brackets. Observe that for $m \geq 2$ 

\[ \frac{q^n (1/q)_n (-1)^m}{q^{rm^2+{m \choose 2}} (1/q)_m} Coeff. \ u^{n-rm} \ in \ (1-u/q^{2m}) (u/q)_{m-1} \prod_{i \geq 1} \frac{1}{1-u/q^i} \] is negative only when $m$ is odd, in which case as above it is less than
\[ \frac{1}{(1-1/q)^2} \frac{1}{q^{{m \choose 2} + n(m-1)}}.\] This gives a contribution of at most

\begin{eqnarray*}
\frac{1}{(1-1/q)^2} \sum_{m \geq 3 \atop m \ odd} \frac{1}{q^{{m \choose 2} + n(m-1)}} & = & \frac{1}{(1-1/q)^2} \sum_{m \geq 1} \frac{1}{q^{2mn+2m^2+m}}\\
& \leq & \frac{1}{(1-1/q)^2} \sum_{m \geq 1} \frac{1}{q^{m(2n+3)}}\\
& \leq & \frac{1}{(1-1/q)^2} \frac{1}{q^{2n+3}(1-1/q^{2n+3})}.
\end{eqnarray*}

\end{proof}
	
	We conclude this section by proving the monotonicity result that $P_{n,q}^r \geq P_{n+1,q}^r$ if $q \geq 2$. Although this result will not be needed elsewhere in the paper, it is combinatorially interesting and may be useful in the future. An analogous result exists for Plancherel measure \cite{Jo} and was crucial for the dePoissonization step in understanding the distribution of the longest increasing subsequence of a random permutation \cite{BDJ},\cite{BOO}. For the case of Plancherel measure, the monotonicty result is true because there is a simple growth process for generating the random partitions such that at stage n of the process has the correct distribution on partitions of size n. Although there is a method for sampling from $P_{n,q}$ (Section 3.3 of \cite{F2}), it is not evident how it can be used to prove the monotonicity result.

	To proceed we require a tool. Recall the Young
Lattice: the elements of this lattice are all partitions of all
natural numbers and an edge is drawn between partitions $\lambda$ and $\Lambda$
if $\Lambda$ is obtained from $\lambda$ by adding one dot. 

\begin{theorem} \label{path} (\cite{F1}) Put weights $m_{\lambda,\Lambda}$ on the Young lattice according to the rules:

\begin{enumerate}

\item $m_{\lambda,\Lambda} = \frac{1}{q^{\lambda_1'}(q^{\lambda_1'+1}-1)}$ if
$\Lambda$ is obtained from $\lambda$ by adding a dot to column 1

\item $m_{\lambda,\Lambda} = \frac{(q^{-\lambda_s'}-q^{-
\lambda_{s-1}'})}{q^{\lambda_1'}-1}$ if $\Lambda$ is obtained from
$\lambda$ by adding a dot to column $s>1$

\end{enumerate} Then

\[ \frac{1}{\prod_j q^{(\lambda_j')^2} (\frac{1}{q})_{m_j(\lambda)}} = \sum_{\gamma} \prod_{i=0}^{|\gamma|-1} m_{\gamma_i,\gamma_{i+1}}\] where $\gamma=\gamma_0 \mapsto \gamma_1 \mapsto \cdots \mapsto \gamma_{n}=|\lambda|$ is a path in the Young lattice from the empty partition to $\lambda$. \end{theorem}

	We remark in passing that because of Theorem \ref{path} the
	measure $P_{n,q}$ can be refined to give a measure on standard
	Young tableaux of size $n$ (which is the same as a path in the
	Young lattice from the empty partition to a partition of size
	$n$). These tableau correspond to involutions in the symmetric
	group via the Robinson-Schensted-Knuth correspondence, and
	there has been much interest in increasing subsequences in
	involutions (e.g. \cite{BaR} and the applications referenced
	there). It remains to be seen whether the measure arising from
	Theorem \ref{path} has similar applications (for a group theoretic
	application, see Section 3.2 of \cite{F2}).

\begin{theorem} \label{monotone} If $q \geq 2$ then $P_{n,q}^r \geq P_{n+1,q}^r$. 
\end{theorem}

\begin{proof} From Proposition \ref{Pform} it is enough to show that

\[ q^{n+1} (1/q)_{n+1} \sum_{|\lambda|=n+1 \atop \lambda_1<r} \frac{1}{\prod_j q^{(\lambda_j')^2} (\frac{1}{q})_{m_j(\lambda)}} \leq q^n (1/q)_n \sum_{|\lambda|=n \atop \lambda_1<r} \frac{1}{\prod_j q^{(\lambda_j')^2} (\frac{1}{q})_{m_j(\lambda)}}.\] By Proposition \ref{path}, 

\[ q^{n+1} (1/q)_{n+1} \sum_{|\lambda|=n+1 \atop \lambda_1<r} \frac{1}{\prod_j q^{\sum (\lambda_j')^2} (\frac{1}{q})_{m_j(\lambda)}} \leq q^{n+1} (1/q)_{n+1} \sum_{|\lambda|=n \atop \lambda_1<r} \frac{1}{\prod_j q^{(\lambda_j')^2} (\frac{1}{q})_{m_j(\lambda)}} \sum_{\Lambda:\lambda \rightarrow \Lambda} m_{\lambda,\Lambda}.\]

	Thus it is enough to show that for all $\lambda$ of size $n$ with $\lambda_1<r$, 
\[ q(1-1/q^{n+1}) \sum_{\Lambda:\lambda \rightarrow \Lambda} m_{\lambda,\Lambda} \leq 1.\] This is visibly true if $\lambda_1'=0$ (i.e. if $\lambda$ is the empty partition). For $|\lambda| \geq 1$, it is easy to see that $\sum_{\Lambda:\lambda \rightarrow \Lambda} m_{\lambda,\Lambda}=\frac{1}{q^{\lambda_1'}}(1+\frac{1}{q^{\lambda_1'+1}-1})$. Thus it must be shown that for all $\lambda$ of size $n$ with $\lambda_1<r$, 
\[ q(1-1/q^{n+1}) \frac{1}{q^{\lambda_1'}}(1+\frac{1}{q^{\lambda_1'+1}-1}) \leq 1.\] If $\lambda_1'=1$ and $r<n+1$ this holds since no such $\lambda$ exist. Similarly if $r>n+1$ the theorem is true since both probabilities are 1. If $\lambda_1'=1$ and $r=n+1$, we have a problem but are saved since the only legal way to add a dot in such a way as to keep $\lambda_1<r$ is to add to column 1 and $q(1-1/q^{n+1}) \frac{1}{q(q^2-1)}  \leq 1$. Finally if $\lambda_1'\geq 2$, the result follows because \[ q(1-1/q^{n+1}) \frac{1}{q^2}(1+\frac{1}{q^{3}-1}) \leq 1\] for $q \geq 2$. \end{proof}

\section{The Measure $Q_{n,q}$ on Partitions} \label{subsequences}

	To begin we recall the measure $Q_{n,q}$ on partitions of size
	$n$ introduced in \cite{F0} and indicate its significance. The
	measure $Q_{n,q}$ arises by any of the following constructions
	and is a natural $q$-analog of the Plancherel measure of the
	symmetric group.

\begin{enumerate}

\item (\cite{F0},\cite{F1}) Choose a partition $\lambda$ of $n$ with probability proportional to the square of the degree of the unipotent representation of $GL(n,q)$ indexed by $\lambda'$. The normalizing constant is

\[ q^{n^2} (1/q)_n^2 Coeff. \ u^n \ in \ \prod_{i=1}^{\infty} \prod_{j=0}^{\infty} (\frac{1}{1-\frac{u}{q^{i+j}}}).\]

\item (\cite{F1}) Recall that the major index of a permutation $\pi \in S_n$ is defined by \[ maj(\pi) = \sum_{i: 1 \leq i \leq n-1 \atop \pi(i)>\pi(i+1)} i.\] Consider the non-uniform measure on the symmetric group which chooses a permutation $\pi$ with probability proportional to $q^{maj(\pi)+maj(\pi^{-1})}$. Let $\lambda$ be the transpose of the partition associated to $\pi$ through the Robinson-Schensted-Knuth (RSK) correspondence. Note that the first row of this $\lambda$ is the length of the longest decreasing subsequence of $\pi$ and that the first column of this $\lambda$ is the length of the longest increasing subsequence of $\pi$. Equivalently, the first rows and columns of $\lambda$ correspond to longest increasing and decreasing subsequences in the reversal of $\pi$. For background on the RSK correspondence including connections with increasing subsequences, see Chapter 7 of \cite{Sta}.

\item There is a measure $\tilde{Q}_q$ on the set of all partitions of natural numbers which chooses a partition $\lambda$ with probability 

\[ \prod_{i=1}^{\infty} \prod_{j=0}^{\infty} (1-\frac{1}{q^{i+j}}) \frac{1}{\prod_j q^{(\lambda_j')^2} \prod_{s \in \lambda} (1-\frac{1}{q^{h(s)}})^2} \] where $h(s)$ is the hooklength of $s$. The measure $Q_{n,q}$ is given by renormalizing $\tilde{Q}_q$ to live on partitions of size $n$. 

\end{enumerate}

	The first two constructions motivate the study of $Q_{n,q}$. The third construction is what will be used in the remainder of this article so we make some remarks about it before continuing.

{\bf Remarks:}

\begin{enumerate}

\item It is easy to see that the measure $\tilde{P}_q$ on the set of all partitions of natural numbers can be rewritten as

\[ \prod_{i=1}^{\infty} (1-\frac{1}{q^{i}}) \frac{1}{q^{|\lambda|+2n(\lambda)} \prod_{s \in \lambda \atop a(s)=0} (1-\frac{1}{q^{h(s)}})}.\] Thus one sees a striking similarity between $\tilde{Q}_q$ and $\tilde{P}_q$ which was one of the motivations for this article.

\item The measure $\tilde{Q}_q$ chooses a partition with probability
proportional to \[s_{\lambda'}(1,\frac{1}{q},\frac{1}{q^2},\cdots)
s_{\lambda'}(\frac{1}{q},\frac{1}{q^2},\frac{1}{q^3},\cdots)\] and hence
was also studied by Okounkov \cite{O}, who computed ``correlation
functions'' for such measures.  As is clear from \cite{G}, the
$\tilde{Q}_q$ probability of having $\lambda_1<r$ can be expressed as
a Toeplitz determinant.

\item It is possible to exactly sample from all four measures $\tilde{P}_q,\tilde{Q}_q,P_{q,n},Q_{q,n}$. See $\cite{F2}$ for discussion. The cases $P_{q,n},Q_{q,n}$ are joint work with Mark Huber.

\end{enumerate}

	Before delving into a study of subsequences, we observe that the measure $Q_{n,q}$ has a curious symmetry property.

\begin{theorem} $Q_{n,q}(\lambda)=Q_{n,1/q}(\lambda')$.
\end{theorem}

\begin{proof} It is proved in \cite{F0} that $Q_{n,q}$ chooses $\lambda$ with probability proportional to $[s_{\lambda'}(1,\frac{1}{q},\frac{1}{q^2},\cdots)]^2$. From the description of $Q_{n,q}$ in terms of major index, one sees from Proposition 7.19.11 and Lemma 7.23.1 of \cite{Sta} that $Q_{n,q}$ picks $\lambda$ with probability proportional to $[s_{\lambda}(1,q,q^2,\cdots)]^2$. \end{proof}

	Proposition \ref{form} gives a formula for $Q_{q,n}$, which is supported on partitions of size $n$.

\begin{prop} \label{form} (\cite{F0}) \[ Q_{q,n}(\lambda) = \frac{1}{Coeff. \ of \ u^n \ in \  \prod_{i=1}^{\infty} \prod_{j=0}^{\infty} (\frac{1}{1-\frac{u}{q^{i+j}}})}
 \frac{1}{q^{|\lambda|+2n(\lambda)} \prod_{s \in \lambda}
 (1-\frac{1}{q^{h(s)}})^2}. \] \end{prop}

	Our next goal (Lemma \ref{boundconst}) is an upper and lower bound on the normalization constant of the measure $Q_{q,n}$. For this a more preliminary lemma is needed.

\begin{lemma} \label{prelim} If $q \geq 2$ then \[ \prod_{i=1}^{\infty} \prod_{j=0}^{\infty} (1-\frac{1}{q^{i+j}}) \geq (1-1/q)^4.\] \end{lemma}

\begin{proof} Since $q>1$ it follows that

\[ \prod_{i=n}^{\infty} (1-\frac{1}{q^i}) \geq  1-\sum_{i=n}^{\infty} \frac{1}{q^i} = 1-\frac{1}{q^n(1-1/q)}.\] Since $q \geq 2$, $1-\frac{1}{q^n(1-1/q)} \geq 1-\frac{1}{q^{n-1}}$. Thus Lemma \ref{neumann1} gives that

\begin{eqnarray*}
\prod_{i=1}^{\infty} \prod_{j=0}^{\infty} (1-\frac{1}{q^{i+j}}) & \geq & (1-1/q)^2 \prod_{i=2}^{\infty} \prod_{j=0}^{\infty} (1-\frac{1}{q^{i+j}})\\
& \geq & (1-\frac{1}{q})^2 \prod_{i=1}^{\infty} (1-\frac{1}{q^i})\\
& \geq & (1-\frac{1}{q})^4.
\end{eqnarray*}
\end{proof}

\begin{lemma} \label{boundconst} Let $z(n,q)$ be denote the coefficient of $u^n$ in $\prod_{i=1}^{\infty} \prod_{j=0}^{\infty} (\frac{1}{1-\frac{u}{q^{i+j}}})$. Then for $q \geq 2$, \[ \frac{1}{q^n (1/q)_n} \leq z(n,q) \leq \frac{1}{(q^n-1)(1-1/q)^6} .\] \end{lemma}

\begin{proof} For the lower bound, observe that

\[ z(n,q) \geq Coeff. \ of \ u^n \ in \ \prod_{i=1}^{\infty} (\frac{1}{1-\frac{u}{q^{i}}}) = \frac{1}{q^n (1/q)_n}.\] For the upper bound, we begin with the recurrence proved in \cite{F0} that

\begin{eqnarray*}\
z(n,q) & = & \frac{1}{q^n-1} \sum_{i=1}^n \frac{z(n-i,q)}{(1/q)_i}\\
& \leq & \frac{1}{(q^n-1)(1-1/q)^2} \sum_{i=1}^n z(n-i,q)\\
& \leq & \frac{1}{(q^n-1)(1-1/q)^2} \sum_{i=0}^{\infty} z(i,q)\\
& = & \frac{1}{(q^n-1)(1-1/q)^2} \prod_{i=1}^{\infty} \prod_{j=0}^{\infty} (\frac{1}{1-\frac{1}{q^{i+j}}}).
\end{eqnarray*} The result now follows from Lemma \ref{prelim}. \end{proof}

	Proposition \ref{compare} gives upper and lower bounds for $Q_{n,q}$ in terms of $P_{n,q}$.

\begin{prop} \label{compare} For $q \geq 2$, \[ (1-1/q^n) (1-1/q)^4 P_{n,q}(\lambda) \leq Q_{n,q}(\lambda) \leq \frac{1}{(1-1/q)^{-1+4 \sqrt{2n}}}  P_{n,q}(\lambda) .\] \end{prop}

\begin{proof} Assume that $n>0$, the case $n=0$ being clear. For the lower bound, Lemma \ref{boundconst} implies that

\begin{eqnarray*}
Q_{n,q}(\lambda) & \geq & \frac{(q^n-1)(1-1/q)^6}{q^n (1/q)_n} \frac{q^n (1/q)_n}{\prod_j q^{(\lambda_j')^2} \prod_{s \in \lambda} (1-1/q^{h(s)})^2}\\
& = & \frac{(q^n-1)(1-1/q)^6}{q^n (1/q)_n} \frac{1}{\prod_{s \in \lambda: a(s)=0} (1-1/q^{h(s)}) \prod_{s \in \lambda: a(s) \neq 0} (1-1/q^{h(s)})^2} P_{n,q}(\lambda)\\
& \geq & \frac{(1-1/q^n) (1-1/q)^5}{(1/q)_n} P_{n,q}(\lambda)\\
& \geq & (1-1/q^n) (1-1/q)^4 P_{n,q}(\lambda).
\end{eqnarray*} The second inequality uses the fact that for non-empty partitions, there is at least one dot satisfying $a(s) \neq 0$. 

	For the upper bound, Lemma \ref{boundconst} implies that
\begin{eqnarray*}
Q_{n,q}(\lambda) & \leq & \frac{q^n (1/q)_n}{\prod_j q^{(\lambda_j')^2} \prod_{s \in \lambda} (1-1/q^{h(s)})^2}\\
& = & \frac{1}{\prod_{s \in \lambda: a(s)=0} (1-1/q^{h(s)}) \prod_{s \in \lambda: a(s) \neq 0} (1-1/q^{h(s)})^2} P_{n,q}(\lambda)\\
& = & \frac{\prod_{s \in \lambda: a(s)=0} (1-1/q^{h(s)})}{\prod_{s \in \lambda} (1-1/q^{h(s)})^2} P_{n,q}(\lambda)\\
& \leq & \frac{(1-1/q)} {\prod_{s \in \lambda} (1-1/q^{h(s)})^2} P_{n,q}(\lambda).
\end{eqnarray*} Since the number of dots $s$ in $\lambda$ with $h(s)=1$ is equal to the number of distinct parts of $\lambda$, it is at most $\sqrt{2n}$. Removing the dots with $h(s)=1$ and applying the same reasoning shows that the number of dots $s$ in $\lambda$ with $h(s)=2$ is at most $\sqrt{2n}$, and that generally the number of dots $s$ in $\lambda$ with any prescribed $h(s)$ value is at most $\sqrt{2n}$. Thus

\[ \frac{(1-1/q)}{\prod_{s \in \lambda} (1-1/q^{h(s)})^2} P_{n,q}(\lambda)  \leq \frac{(1-1/q)}{\prod_{i} (1-1/q^i)^{2 \sqrt{2n}}} \leq  \frac{1}{(1-1/q)^{-1+4 \sqrt{2n}}}.\]
\end{proof}

	Combining Theorem \ref{Pbound} and Proposition \ref{compare}, we obtain the following result.

\begin{theorem} Let $Q_{n,q}^r$ be the $Q_{n,q}$ probability that $\lambda_1<r$. For $q \geq 2$ and $r \leq n-1$,

\begin{enumerate}
\item $Q_{n,q}^r \leq  \frac{1}{(1-1/q)^{1+4 \sqrt{2n}}} [\frac{1}{q^{2n-2r+2}}+ \frac{1}{(1-1/q^{2n+1})}\frac{1}{q^{n+1}}]$

\item $Q_{n,q}^r \geq (1-1/q^n) (1-1/q)^4 [\frac{1}{q^{2n-2r+2}}-\frac{1}{q^{n+1}(1-1/q)} - \frac{1}{q^{2n-2r+3}(1-1/q)^2} - \frac{1}{q^{3n-3r+4}(1-1/q)^3} - \frac{1}{(1-1/q)^2} \frac{1}{q^{2n+3}(1-1/q^{2n+3})}]$

\end{enumerate}
\end{theorem}

	From Theorem \ref{firstcolumn} and Proposition \ref{compare} one deduces the following corollary.

\begin{cor} (\cite{F1}) Suppose that $q \geq 2$.

\begin{enumerate} 

\item The $Q_{n,q}$ probability that $\lambda_1'=k$ is at most $ \frac{1}{(1-1/q)^{1+4 \sqrt{2n}}} \frac{(1/q)_n (1/q)_{n-1}}{q^{k^2-k} (1/q)_k (1/q)_{k-1} (1/q)_{n-k}}$.

\item The $Q_{n,q}$ probability that $\lambda_1'=k$ is at least $(1-1/q^n) (1-1/q)^4 \frac{(1/q)_n (1/q)_{n-1}}{q^{k^2-k} (1/q)_k (1/q)_{k-1} (1/q)_{n-k}}$.

\end{enumerate}
\end{cor}


\begin{thebibliography}{AAA}

\bibitem [AD] {AD} Aldous, D. and Diaconis, P., Longest increasing subsequences: from patience sorting to the Baik-Deift-Johansson theorem, Bull. AMS (N.S.) {\bf 36} (1999), 413-432.

\bibitem [BaDeJo] {BDJ} Baik, J., Deift, P., and Johansson, K., On the length of the longest increasing subsequence of random permutations, {\it J. Amer. Math. Soc.} {\bf 12} (1999), 1119-1178.

\bibitem [BaR] {BaR} Baik, J. and Rains, E., The asymptotics of monotone subsequences of involutions. Available at http://xxx.lanl.gov/abs/math.CO/9905084.

\bibitem [BOOl] {BOO} Borodin, A., Okounkov, A., and Olshanski, G., Asymptotics of Plancherel measures for symmetric groups, {\it J. Amer. Math. Soc.} {\bf 13} (2000), 481-515.

\bibitem [De]{De} Deift, P., Integrable systems and combinatorial theory, {\it Notices Amer. Math. Soc.} {\bf 47} (2000), 631-640.

\bibitem [F0] {F0} Fulman, J., {\it Probability in the classical groups over finite fields: symmetric functions, stochastic algorithms and cycle indices}, Ph.D. Thesis, Harvard University, 1997.

\bibitem [F1] {F1} Fulman, J., A probabilistic approach to conjugacy classes in the finite general linear and unitary groups, {\it J. Algebra} {\bf 212} (1999), 557-590.

\bibitem [F2] {F2} Fulman, J., Random matrix theory over finite fields. To appear in {\it Bull. AMS (N.S.)}. Available at http://www.math.pitt.edu/$\sim$fulman.

\bibitem [F3] {F3} Fulman, J., A probabilistic proof of the Rogers-Ramanujan identities, {\it Bull. London Math. Soc} {\bf 33} (2001), 397-407.

\bibitem [G] {G} Gessel, I., Symmetric functions and P-recursiveness, {\it J. Combin. Theory Ser. A} {\bf 53} (1990), 257-285.

\bibitem [H] {H} Herstein, I.N., Topics in algebra, 2nd edition. Xerox
Corporation. 1975.

\bibitem [Jo] {Jo} Johansson, K., The longest increasing subsequence in a random permutation and a unitary random matrix model, {\it Math. Res. Lett.} {\bf 5} (1998), 63-82. 

\bibitem [NP] {NP} Neumann, P.M. and Praeger. C.E., Cyclic matrices over finite fields, {\it J. London Math. Soc. (2)} {\bf 52} (1995), 263-284.

\bibitem [O]{O} Okounkov, A., Infinite wedge and random partitions. Available at http://xxx.lanl.gov/abs/math.RT/9907127.

\bibitem [S]{Sta} Stanley, R., Enumerative combinatorics (Volume 2). Cambridge University Press, Cambridge, UK. 1999.

\bibitem [TW]{TW} Tracy, C. and Widom, H., On the distributions of the lengths of the longest monotone subsequences in random words, {\it Probab. Theory Related Fields} {\bf 119} (2001), 350-380.
		
\end{thebibliography}
\end{document}